\input amstex
\documentstyle{amsppt}
\NoBlackBoxes
\pagewidth{6.4truein}
\pageheight{9truein}
\def\cH{{\Cal H}}

\def\complex{{\Bbb C}}
\def\IH{{\Bbb H}}
\def\nat{{\Bbb N}}
\def\real{{\Bbb R}}

\def\ep{\varepsilon}
\def\Bone{\text{\bf 1}}
\topmatter
\title On the Grushin operator and hyperbolic symmetry 
\endtitle
\author William Beckner \endauthor 
\address Department of Mathematics, The University of Texas at Austin, 
Austin, TX 78712-1082\endaddress 
\email beckner\@math.utexas.edu\endemail
\thanks This work was partially supported by the National Science 
Foundation.\endthanks 
\abstract
Complexity of geometric symmetry for differential operators with
mixed homogeniety is examined here. Sharp Sobolev estimates are 
calculated for the Grushin operator in low dimensions using 
hyperbolic symmetry and conformal geometry.
\endabstract 
\endtopmatter
  
\document 

Considerable interest exists in understanding differential operators with 
mixed homogeneity. 
A simple example is the Grushin operator on $\real^2$ 
$$\Delta_G = {\partial^2\over\partial t^2} 
+ 4t^2 {\partial^2\over\partial x^2}\ .$$ 
The purpose of this note is to demonstrate the complexity of geometric 
symmetry that may exist for operators defined on Lie groups. 
Here the existence of an underlying $SL(2,R)$ symmetry for $\Delta_G$ is 
used to compute the sharp constant for the associated $L^2$ Sobolev inequality. 

\proclaim{Theorem 1} 
For $f\in C^1 (\real^2)$ 
$$\Big[ \|f\|_{L^6 (\real^2)} \Big]^2 
\le  \pi^{-2/3} \int_{\real^2} 
\left[ \Big( {\partial f\over\partial t}\Big)^2 
+ 4t^2 \Big( {\partial f\over\partial x}\Big)^2\right]\, dx\,dt\ .
\tag 1$$ 
This inequality is sharp, and an extremal is given by
$\big[ (1+|t|^2)^2+|x|^2\big]^{-1/4}$.
\endproclaim 

This result follows from the analysis of a Sobolev inequality on 
$SL(2,R)/ SO(2)$. But the hyperbolic embedding estimate requires some
interpretation to take into account cancellation effects. It will be
essential to include contibutions to the hyperbolic Dirichlet form from
non-$L^2$ functions.
Let $z= x+iy$ denote a point in the upper half-plane $\real_+^2\simeq 
\IH^2 \simeq M\simeq SL(2,R)/SO(2)$. 
Here the invariant distance is given by the Poincar\'e metric 
$$d(z,z') = {|z-z'| \over 2\sqrt{yy'}}$$ 
with the corresponding invariant gradient $D=y\nabla$ and left-invariant 
Haar measure $d\nu = y^{-2} \,dy\,dx$. 

\proclaim{Theorem 2}
For $F\in C^1_c(M)$ 
$$\Big[ \|F\|_{L^6(M)}\Big]^2 \le 4 \pi^{-2/3}
\biggl[ \int_M |DF|^2\,d\nu -\dfrac3{16} \int_M |F|^2\,d\nu\biggr]\  
\tag 2$$
$$\Big[ \|F\|_{L^6(M)}\Big]^2 \le \dfrac43 \pi^{-2/3}
\biggl[ \int_M |DF|^2\,d\nu - \dfrac14 \int_M |F|^2\,d\nu\biggr]\ .
\tag 3$$ 
Both estimates are sharp as limiting forms.
\endproclaim 

These two estimates would seem to be contradictory,
but it must be understood that the right-hand sides are to be evaluated 
as limiting forms for functions
that may not be in $L^2(M)$. So the issue of which is the sharper Sobolev
inequality must be studied carefully. On the hyperbolic manifold the Dirichlet
form can be represented as a weighted Sobolev form so that for $\alpha > 0$
$$
\int_M y^{2\alpha} |\nabla y^{-\alpha} f|^2\,dx\,dy 
= \int_M |Df|^2\, d\nu + \alpha (\alpha - 1) \int_M |f|^2\,d\nu\ .
$$
On the right-hand side of equation (3) observe the appearance of the spectral
limit $\frac 14$ for the hyperbolic Laplacian $- y^2 \Delta$:
$$\dfrac14 \int_M |F|^2\,d\nu \le \int_M |DF|^2\,d\nu\ .
\tag 4$$

\demo{Proof of Theorem 1} 
Let $\tilde f$ denote the Fourier transform of $f$ in the first variable. 
That is, for integrable functions 
$$\tilde f(\xi,t) = \int_{\real} e^{2\pi i\xi x} f(x,t)\,dx$$ 
so that by using the Plancherel identity, inequality (1) for some constant 
$A_0$ is equivalent to 
$$\biggl[ \int_{\real^2} |\tilde f*\tilde f*\tilde f|^2\,d\xi\, dt\bigg]^{1/3}
\le A_0 \int_{\real^2} \biggl[ \Big( {\partial \tilde f\over \partial t}
\Big)^2 + 16\pi^2 |t|^2 |\xi|^2 |\tilde f|^2 \biggr]\,d\xi\,dt$$ 
where here convolution is only with respect to the first variable. 
Now one can apply standard rearrangement arguments of Riesz-Sobolev type 
to see that it suffices to consider this inequality only for non-negative 
functions $\tilde f(\xi,t)$ that are symmetric decreasing in each of the 
two variables separately. 
Hence, the function $f(x,t)$ in (1) can be taken to be symmetric in $t$ and 
symmetric decreasing in $x$. 
The second part of this remark follows from the fact that the Dirichlet 
form in (1) taken only with respect to integration in $x$ is diminished 
by a symmetric decreasing equimeasurable rearrangement in the first variable. 

Since $f(x,t)$ is even in $t$, set $y=t^2$ and 
let $f(x,|t|) = y^{-1/4}F(x,y)$; then 
$$\|f\|_{L^6(\real^2)} = \|F\|_{L^6(M)}$$ 
and inequality (1) is now equivalent to 
$$\Big[ \|F\|_{L^6(M)} \Big]^2 \le 4A_0 
\biggl[ \int_M |DF|^2 \,d\nu - \dfrac3{16} \int_M |F|^2\,d\nu\biggr]\ .$$ 
This is an {\it a~priori\/} 
inequality where the function $F$ can be taken to be 
smooth but still the form will extend to functions that are not in $L^2(M)$.
One can also restrict this result 
to consideration of functions that are radial decreasing in the Poincar\'e 
distance from the origin $\hat 0 = (0,1)=i$. 
Now Theorem~1 will follow from the first part of Theorem~2 with 
$A_0 = \pi^{-2/3}$.
\enddemo 

\demo{Proof of equation (2) in Theorem 2}
By using equimeasurable radial decreasing rearrangement corresponding to the 
metric on hyperbolic space, it suffices to consider this inequality for 
radial decreasing functions of the distance from the origin. 
Let $u= [d(z,i)]^2$; then for functions depending on distance the 
gradient is given by 
$$|DF| = \sqrt{u+u^2}\, \Big| {dF \over du}\Big|$$ 
and the volume form restricted to integration for radial function is given 
by $d\nu = 4\pi\,du$. 
Then (2) is equivalent to 
$$\biggl[ \int_0^\infty |F|^6\,du\bigg]^{1/3} 
\le 2^{10/3} \biggl[  \int_0^\infty (u+u^2) |F'|^2\,du 
- \dfrac3{16} \int_0^\infty |F|^2\,du\biggr]\ .$$ 
Let $G\in C_c^2 ([0,\infty))$ and set $F(u) = (1+u)^{-1/4} G(u)$. 
Then inequality (2) is equivalent to 
$$\biggl[ \int_0^\infty |G|^6 (1+u)^{-3/2}\,du\biggr]^{1/3} 
\le 2^{10/3} \bigg[ \int_0^\infty u\sqrt{1+u}\, |G'|^2\,du 
+ \dfrac1{16} \int_0^\infty |G|^2 (1+u)^{-3/2}\,du\biggr]\ .$$
But now this estimate will be considered for all Lipschitz
functions $G$ such that the right-hand side is finite.
Make the change of variables $u\to 1/u$ with $H(u) = G(1/u)$; 
$$\eqalign{
\biggl[ \int_0^\infty |H|^6 (1+u)^{-3/2} u^{-1/2}\,du\biggr]^{1/3} 
&\le 2^{10/3} \bigg[ \int_0^\infty \sqrt{u(1+u)}\, |H'|^2 \,du \cr 
&\qquad + \dfrac1{16} \int_0^\infty |H|^2 (1+u)^{-3/2} u^{-1/2}\,du\biggr]\ .
\cr}\tag 5$$
By evaluating this estimate for $H(u) = (1+u)^{-\ep}$ as $\ep\to0$, 
one sees that the constant cannot be smaller than $2^{10/3}$. 
This calculation also suggests that the inequality should be associated 
with sharp Sobolev embedding on $S^2$. Such intuition is realized by the
following argument.

Define a new variable $w$ by setting 
$$(1+w)^{-2} \,dw = \dfrac12 u^{-1/2} (1+u)^{-3/2}\,du$$ 
so that 
$$\sqrt{ u\over 1+u} = {w\over 1+w}$$ 
and $w= u+\sqrt{u(1+u)}$. 
With this change of variables (5) becomes 
$$\biggl[ \int_0^\infty |G|^6 (1+w)^{-2} \,dw\biggr]^{1/3} 
\le 4 \int_0^\infty (2w+1) |G'|^2\,dw 
+ \int_0^\infty |G|^2 (1+w)^{-2} \,dw\ .
\tag 6$$ 
This inequality is controlled by sharp Sobolev embedding on $S^2$; 
more precisely, the family of sharp Sobolev inequalities on $S^2$ that are 
determined by the Hardy-Littlewood-Sobolev inequality (see Theorem~4 in 
\cite2) 
$$\biggl[ \int_{S^2} |F|^p\,d\xi\biggr]^{2/p} 
\le {p-2\over2} \int_{S^2} |\nabla F|^2\,d\xi 
+ \int_{S^2} |F|^2\,d\xi 
\tag 7$$ 
for $2\le p<\infty$ and $d\xi$ is normalized surface measure on $S^2$. 
Inequality (6) follows from the case $p=6$. 
Observe that the change of variables defined by stereographic projection 
between $\real^2$ and $S^2$-\{pole\} can be realized for the polar angle 
on $S^2$ by $\cos\theta = (1-|x|^2)/(1+|x|^2)$ and $w= |x|^2$ in (6). 
Since inequality (6) then corresponds to functions of the polar angle, it 
suffices simply to match up the ``radial coordinates'' in each domain. 
Then 
$$w\Big( {dG\over dw}\Big)^2 \,dw = \dfrac12 \Big( {dG\over d\theta}\Big)^2
\sin\theta\,d\theta$$ 
so that (7) for $p=6$ and radial variables gives a stronger inequality 
than (6) 
$$\biggl[ \int_0^\infty |G|^6 (1+w)^{-2} \,dw\biggr]^{1/2} 
\le 2 \int_0^\infty w |G'|^2\,dw 
+ \int_0^\infty  |G|^2 (1+w)^{-2}\,dw\ .
\tag 8$$
This shows that inequality (6) is sharp only for constants.  
\enddemo

\demo{Proof of equation (3) in Theorem 2} 
This result is a special case of an argument in \cite4 that uses axial 
symmetry and $SL(2,R)$ to derive the sharp Sobolev embedding constant 
on $\real^n$ and characterize the extremals for that problem. 
The motivation for this approach came from problems in fluid mechanics 
and vortex dynamics. 
For $n>2$ and $1/p = 1/2 - 1/n$ 
$$\gather 
\|f\|_{L^p(\real^n)} \le A_p \|\nabla f\|_{L^2(\real^n)} \tag 9\cr 
A_p = [ \pi n(n-2)]^{-1/2} [\Gamma (n)/\Gamma (n/2)]^{1/n}
\endgather$$ 
and up to the action of the conformal group, the sharp constant is only 
attained for functions of the form $A(1+|x|^2)^{-n/p}$. 
By using the technique of symmetrization (equimeasurable radial decreasing 
rearrangement), it suffices to consider this inequality for non-negative 
radial decreasing functions. 
For radial functions use the product structure for Euclidean space 
$\real^n \simeq \real\times \real^{n-1}$ with $x= (t,x')$ and set $y=|x'|$.  
Being radial in $x$ means that the function is also radial in $x'$. 
Let $g(t,y) = y^{n/p} f(t,x')$ and inequality (9) becomes 
$$\biggl[ \int_M |g|^p\,d\nu\bigg]^{2/p} 
\le B_p \biggl[ \int_M |Dg|^2\,d\nu + \dfrac{n}p \Big(\dfrac{n}p-1\Big) 
\int_M |g|^2\,d\nu\biggr]
\tag 10$$
where 
$$B_p = {4\over n(n-2)} \Big[ {n-1\over 2\pi}\Big]^{2/n}\ .$$ 
For the case $n=3$, $p=6$ and $B_6 = (4/3) \pi^{-2/3}$, and equation (3)
is proved. 
The argument in \cite4 to obtain extremals for the Sobolev inequality (9) is 
a nice application of the competing radial and cylindrical symmetry. 
\enddemo

This result on sharp Grushin estimates is interesting because (1) the 
solution does match the pattern suggested by the Heisenberg group (see 
\cite3), (2) the analysis is controlled by the two-dimensional sharp 
Hardy-Littlewood-Sobolev inequality, and
(3) the identification of $SL(2,R)$ symmetry is related to the role of 
analyticity in the Lewy example. 

For higher dimensions this problem has corresponding behavior. 
Consider $(x,t) \in \real\times \real^2 \simeq \real^3$ with 
$$\Delta_G = \Delta_t + 4|t|^2 {\partial^2\over\partial x^2}\ .$$ 
The homogeneous dimension of this operator is 4. 
Here one can also use the underlying $SL(2,R)$ symmetry to compute the 
sharp constant for the associated $L^2$ Sobolev inequality with a
similar analysis.

\proclaim{Theorem 3} 
For $f\in C^1(\real^3)$ 
$$\Big[ \|f\|_{L^4(\real^3)}\Big]^2 \le {1\over 2\pi} 
\int_{\real\times\real^2} \left[ |\nabla_t f|^2 + 4|t|^2 
\Big( {\partial f\over\partial x}\Big)^2\right]\, dx\,dt\ .
\tag 11$$ 
This inequality is sharp, and an extremal is given by $\big[ (1+|t|^2)^2+
|x|^2\big]^{-1/2}$.
\endproclaim 

\proclaim{Theorem 4}
For $F\in C_c^1(M)$ 
$$\Big[ \|F\|_{L^4 (M)} \Big]^2 \le {2\over\sqrt\pi} 
\biggl[ \int_M |DF|^2\,d\nu  - {1\over4} \int_M |F|^2\,d\nu\biggr]\ . 
\tag 12$$ 
\endproclaim 

\demo{Proof of Theorem 3} 
Let $\tilde f$ denote the Fourier transform of $f$ in the first variable $x$. 
Using the Plancherel identity, inequality (11) for some constant $A_0$ 
is equivalent to 
$$\biggl[\int_{\real^3} |\tilde f* \tilde f|^2 \,d\xi\,dt\biggr]^{1/2} 
\le A_0 \int_{\real^3} \left[ |\nabla_t \tilde f|^2 + 16\pi^2 |t|^2 
|\xi|^2 |\tilde f|^2\right] \,d\xi\,dt$$ 
where here convolution is only with respect to the first variable. 
By applying Riesz-Sobolev rearrangement arguments, it suffices to consider 
this inequality only for non-negative functions $\tilde f(\xi,t)$ that 
are radial decreasing in each of the two variables separately. 
Hence, the function $f(x,t)$ in (11) can be taken to be radial in $t$ 
and symmetric decreasing in $x$. 
The second part of this remark follows from the fact that the Dirichlet 
form in (11) taken only with respect to integration in $x$ is diminished 
by a symmetric decreasing equimeasurable rearrangement with respect to the 
first variable. 

Since $f(x,t)$ is radial in $t$, set $y= |t|^2$ and let $f(x,|t|)=y^{-1/2} 
F(x,y)$; then 
$$\|f\|_{L^4 (\real^3)} = \pi^{1/4} \|F\|_{L^4(M)}$$ 
and inequality (11) is now equivalent to 
$$\Big[ \|F\|_{L^4(M)}\Big]^2 \le 4\sqrt{\pi}\, A_0 
\biggl[ \int_M |DF|^2\, d\nu - {1\over4} \int_M |F|^2\,d\nu\biggr]\ .$$ 
This is an {\it a priori\/} inequality where  the function $F$ can be taken 
to be smooth with compact support. 
Now Theorem~3 will follow from Theorem 4 with $A_0 = 1/(2\pi)$. 
One simply calculates that equality is attained for the indicated extremal. 
\enddemo 

\demo{Proof of Theorem 4}
Using equimeasurable radial decreasing rearrangement corresponding to the 
metric on hyperbolic space, it suffices to consider this inequality for 
radial decreasing functions of the distance from the origin. 
Set $u= d^2 (z,i)$; then the volume form restricted to integration 
for radial functions  
is given by $d\nu = 4\pi \,du$ and inequality (12) becomes (see \cite5) 
$$\biggl[ \int_0^\infty |F|^4\,du\biggr]^{1/2} 
\le 4\biggl[ \int_0^\infty (u^2+u) \Big| {dF\over du}\Big|^2\,du - {1\over4} 
\int_0^\infty |F|^2\,du\biggr]\ .
\tag 13$$ 
If one can show that this is a good upper bound, then the sequence of 
functions $F_\ep(u) = (1+u)^{-\ep}$ for $\ep>\frac12$ shows that the 
estimate is sharp. 
Let $G\in C_c^2 ([0,\infty))$ and set $F(u) = (1+u)^{-1/2} G(u)$. 
Then inequality (13) takes the form 
$$\bigg[ \int_0^\infty |G|^4 {1\over (1+u)^2}\,du\biggr]^{1/2} 
\le 4\int_0^\infty u(G')^2\,du + \int_0^\infty |G|^2 
{1\over (1+u)^2}\, du\ . 
\tag 14$$ 
This inequality is controlled by sharp Sobolev embedding on $S^2$; 
more precisely, the family of sharp Sobolev inequalities on $S^2$ that are 
determined by the Hardy-Littlewood-Sobolev inequality (see Theorem~4 
in \cite2) 
$$\biggl[\int_{S^2} |F|^p\,d\xi\biggr]^{2/p} 
\le {p-2\over2} \int_{S^2} |\nabla F|^2\,d\xi 
+ \int_{S^2} |F|^2\,d\xi
\tag 15$$ 
for $2\le p<\infty$ and $d\xi$ is normalized surface measure on $S^2$. 
Inequality (14) follows from the case $p=4$. 
Observe that change of variables defined by stereographic projection 
between $\real^2$ and $S^2$-$\{$pole$\}$ can be realized for the polar  
angle on $S^2$ by $\cos\theta = {(1-|x|^2)/(1+|x|^2)}$ and $u= |x|^2$ in (14). 
Since inequality (14) corresponds to functions of the polar angle, it suffices 
simply to match up the ``radial coordinates'' in each domain. 
Then 
$$u\Big( {dG\over du}\Big)^2\,du = \frac12 \Big( {dG\over d\theta}\Big)^2 
\sin \theta\,d\theta$$ 
so that (15) for $p=4$ and radial variables gives a stronger 
inequality than (14) 
$$\biggl[ \int_0^\infty |G|^4 {1\over (1+u)^2}\,du\biggr]^{1/2} 
\le \int_0^\infty u(G')^2\,du + \int_0^\infty |G|^2 {1\over (1+u)^2}\,du\ .
\tag 16$$ 
This shows that Theorem 4 is sharp as a limiting form. However, the limit
``extremal'' 
$$F(x,y) = [1 + d^2(z,i)]^{-1/2}$$
is not in $L^2(M)$. This observation emphasizes that the appropriate 
Dirichlet form for Sobolev embedding on hyperbolic space $\Bbb H^2$
should correspond to the intrinsic
positive elliptic differential operator 
$$L_s = -y^2\Delta + s(s-1)\Bone.$$
\enddemo 

These two results illustrate the complexity and interdependence of Sobolev 
estimates on Lie groups and symmetric spaces, and demonstrate that there 
is still much to understand about the geometry of Grushin operators. 
The elementary nature of these calculations was facilitated by the capability 
to use rearrangement arguments which here depended on the Sobolev index 
being an even integer. 
An interesting aspect of the analysis is that the intermediate estimate 
on hyperbolic space must be defined as a limiting form using the positive 
elliptic operator $L_s$ at the extremal for the Grushin embedding estimate. 

\head Appendix\endhead

The argument used here to relate sharp Sobolev embedding on $S^2$ to
embedding estimates on hyperbolic space determines a more general family
of such estimates.

\proclaim {Theorem 5} 
For $F\in C_c^1 (M)$, $0<s\le \frac12$ and $p= 2+\frac1s \ge 4$ 
$$\gather 
\Big[\|F\|_{L^p(M)}\Big]^2 \le A_p \biggl[ \int_M |DF|^2\,d\nu 
+ s(s-1)\int_M |F|^2\,d\nu\biggr] \tag 17\cr
A_p = (2\pi)^{\frac2p -1} s^{-1-\frac2p}\ . 
\endgather$$
\endproclaim

\demo{Proof} 
By using equimeasurable radial decreasing rearrangement corrresponding 
to the metric on hyperbolic space, it suffices to consider this inequality 
for radial decreasing functions of the distance from the origin. 
Let $u=  [d(z,i)]^2$; then (17) is equivalent to 
$$\biggl[ \int_0^\infty |g|^p\,du \biggr]^{2/p} 
\le (4\pi)^{1-\frac2p} A_p \biggl[\int_0^\infty (u+u^2) |g'|^2 \,du 
+ s (s-1) \int_0^\infty |g|^2\,du \biggr]\ .$$
Set $g= (1+u)^{-\alpha} h$ and $C_p = (4\pi)^{(p-2)/p} A_p$; then 
$$\align 
\biggl[ \int_0^\infty |h|^p (1+u)^{-p\alpha}\,du \biggr]^{2/p} 
& \le C_p \biggl[ \int_0^\infty u(1+u)^{1-2\alpha} |h'|^2\,du 
+ \alpha^2 \int_0^\infty |h|^2 (1+u)^{-1-2\alpha}\,du \cr 
&\qquad + (s^2 -s+\alpha - \alpha^2) \int_0^\infty |h|^2 (1+u)^{-2\alpha}\,du 
\biggr]\ .
\endalign$$ 
Set $\alpha =s$, $p\alpha = 2\alpha +1$ and $\beta = 2\alpha >0$; then 
$$\align 
\biggl[ \int_0^\infty |h|^p (1+u)^{-\beta -1}\,du \biggr]^{2/p} 
&\le C_p \biggl[ \int_0^\infty u(1+u)^{1-\beta} |h'|^2\,du \cr
&\qquad + \frac{\beta^2}4 \int_0^\infty |h|^2 (1+u)^{-1-\beta}\,du \biggr]\ .
\endalign $$
Make the change of variables $u\to 1/u$ with $H(u) = h(1/u)$ so that 
$|h'(1/u)| = u^2 |H'(u)|$ and  
$$\align 
\biggl[\int_0^\infty |H|^p (1+u)^{-\beta-1} u^{\beta-1}\,du\biggr]^{2/p} 
&\le C_p \biggl[ \int_0^\infty u^\beta (1+u)^{1-\beta} |H'|^2\,du \cr
&\qquad + \frac{\beta^2}4 \int_0^\infty  |H|^2(1+u)^{-\beta-1} u^{\beta-1}\,du 
\biggr]\ .
\endalign$$
Now set $(1+w)^{-2}\,dw = \beta (1+u)^{-\beta-1} u^{\beta-1}\,du$ so that  
$$\frac{w}{1+w} = \left( \frac{u}{1+u}\right)^\beta$$ 
which gives for $G(w) = H(u)$, $2/p = \beta/(1+\beta)$  and 
$B_p = \frac14 \beta^{(1+2\beta)/(1+\beta)} C_p$ 
$$\displaylines{
(18)
\qquad \biggl[ \int_0^\infty\mkern-9mu 
|G|^p(1+w)^{-2}\,dw \biggr]^{\beta/(1+\beta)}\hfill\cr
\hfill\qquad\qquad \le B_p 
\biggl[ 4 \int_0^\infty \mkern-9mu w^{2-\frac1{\beta}} 
\left[ (1+w)^{1/\beta} -w^{1/\beta} \right] |G'|^2\,dw 
 + \int_0^\infty \mkern-9mu |G|^2 (1+w)^{-2}\,dw \biggr]\ .\hfill}
$$
This equation can be simplified using the change of variables $w\to 1/w$ 
and setting $\widetilde G (w) = G(1/w)$:  
$$\displaylines{
(19)
\qquad \biggl[ \int_0^\infty \mkern-9mu |\widetilde G|^p 
(1+w)^{-2}\,dw\biggr]^{\beta/(1+\beta)}\hfill \cr 
\hfill\qquad\qquad \le B_p \biggl[4\int_0^\infty \left[(1+w)^{1/\beta}-1\right] 
|\widetilde G'|^2\,dw  
+ \int_0^\infty\mkern-9mu  |\widetilde G|^2 (1+w)^{-2}\,dw \biggr]\ .\hfill}
$$
Now this estimate should be compared with the sharp Sobolev embedding on 
$S^2$ that is determined by the Hardy-Littlewood-Sobolev inequality: 
$$\biggl[ \int_{S^2} |F|^p \,d\xi\biggr]^{2/p} 
\le \frac{p-2}2 \int_{S^2} |\nabla F|^2\,d\xi 
+ \int_{S^2} |F|^2 \,d\xi$$ 
where $d\xi$ denotes normalized surface measure and $p>2$; 
and in turn gives for 
radial functions and $p=2 (1+1/\beta)$ 
$$\bigg[ \int_0^\infty |F|^p (1+w)^{-2}\,dw \biggr]^{\beta /(1+\beta)} 
\le \frac1{\beta} \int_0^\infty w|F'|^2\,dw 
+ \int_0^\infty |F|^2 (1+w)^{-2}\,dw \ .
\tag 20$$
Now set $B_p=1$ in (19) which corresponds to the value of $A_p$ in (17) 
and observe that for $r=1/\beta$ and $w\ge0$, 
then $(1+w)^r \ge 1+rw$ for $r\ge 1$. 
Hence the estimate (20) derived from Sobolev embedding on $S^2$ implies 
that (19) holds for $B_p=1$ and $\beta\ge 1$. 
The proof of Theorem~5 is then complete for $0<s\le 1/2$.
\enddemo

In the analysis of Sobolev embedding on the Heisenberg group $\cH_n$ 
realized as 
the manifold  $\complex^n\times\real$ and restricted to radial symmetry in 
the complex variables, then a discrete set of hyperbolic embedding estimates 
can be obtained (see Theorem~18 in \cite3). 

\proclaim{Theorem 6}
For $F\in C^1(M)\cap L^2(M)$, $s=n/2$ for $n\in\nat$ and $p=2+\frac1s \le 4$
$$\gather 
\Big[ \|F\|_{L^p(M)}\Big]^2 
\le A_p\biggl[ \int_M |DF|^2\, d\nu + s(s-1)\int_M |F|^2\,d\nu\cr 
A_p = (2\pi)^{\frac2p -1} s^{-1-\frac2p}\ .
\endgather$$ 
For $n>1$ and up to the ``conformal structure'' of $M$, an extremal 
is given by 
$$F(z) = \Big[ 1+d^2 (z,i)\Big]^{-s}\ .$$
\endproclaim

This family of hyperbolic embedding estimates can be extended to include 
values of $s\ge 1$ by using duality and the fundamental solution 
corresponding to the differential operator $L_s$. 
Note that in this case an $L^2$ extremal function will exist. 
The fundamental solution for $L_s = -y^2\Delta + s(s - 1) 
\Bone$ for $s \ge 1$ is given by 
$$\align 
\psi_s (u) & = \frac1{4\pi} \int_0^1 \big[ t(1-t)\big]^{s-1} (t+u)^{-s}\,dt\cr
&= \frac{\Gamma (s)\Gamma (s)}{4\pi \Gamma (2s)} (1+u)^{-s} 
F\left( s,s,2s;\frac1{1+u}\right)
\endalign$$ 
where $u= [d(z,i)]^2$ and $F$ is the hypergeometric function. 
The transition from Sobolev embedding estimates to a Hardy-Littlewood-Sobolev 
convolution inequality is made using the following lemma. 

\proclaim{Lemma} 
Let $K$ and $\Lambda$ be densely defined, positive-definite, 
self-adjoint operators acting 
on functions defined on a $\sigma$-finite measure space $M$ and satisfying 
the relation 
$$\Lambda K = K\Lambda = \text{\bf 1}\ .$$
Then the following two inequalities are equivalent: 
$$\gather 
\|Kf \|_{L^{p'}(M)} \le C_p \|f\|_{L^p(M)} \tag $*$\cr 
\|g\|_{L^{p'}(M)} \le \sqrt{C_p} \|\Lambda^{1/2} g\|_{L^2(M)}\ .\tag $**$
\endgather$$
Here $1<p<2$ and $1/p + 1/p' =1$. 
Extremal functions for one inequality will determine extremal functions 
for the other inequality if the operator forms are well-defined. 
\endproclaim

\demo{Proof} 
In $(**)$ substitute $g=Kf$ so that 
$$\align 
\Big[ \|Kf\|_{L^{p'}(M)}\Big]^2 
&\le C_p \int_M (Kf) \Lambda (Kf) \,dm = C_p \int_M (Kf)f\,dm \\
&\le C_p \|Kf\|_{L^{p'}(M)} \|f\|_{L^p(M)}
\endalign$$
which is now $(*)$. 
For equivalence in the reverse direction, $K$ is a positive-definite 
self-adjoint operator and notice that $(*)$ implies 
$$\|K^{1/2} f\|_{L^2 (M)} \le \sqrt{C_p} \|f\|_{L^p(M)}$$
which by duality implies 
$$\|K^{1/2} h\|_{L^{p'}(M)} \le \sqrt{C_p} \|h\|_{L^2 (M)}\ .$$
Now substitute $h=K^{1/2} (\Lambda g)$ which results in $(**)$. 
The full equivalence is obtained by taking limits on dense domains. 

For $s>0$ define the fractional integral operator 
$$(I_sG) (z) = \int_M \psi_s \big[ d^2 (z,w)\big] G(w)\,d\nu\ . 
\tag 21$$
The symmetric space $SL(2,R)/SO(2)\simeq \IH^2$ can be identified 
with the subgroup of $SL (2,R)$ given by all matrices of the form 
$$\pmatrix \sqrt y & x/\sqrt y\\ 
\noalign{\vskip6pt}
0&1/\sqrt y 
\endpmatrix$$
with $y>0$ and $x\in\real$ which act via fractional linear transformations 
on $\real_+^2 \simeq \IH^2$. 
$$z = x+iy \in\real_+^2 \to \frac{az+b}{cz+d}$$ 
for $\left(\smallmatrix a&b\\ c&d\endsmallmatrix\right) \in SL(2,R)$. 
The modular function is $\Delta (x,y) = 1/y$ and $d\nu = y^{-2}\,dx\,dy$ 
is left-invariant Haar measure on the group. 
Observe that the group action here corresponds to the multiplication rule 
$$(x,y) (u,v) = (x+yu,yv)$$ 
for $x,u\in\real$ and $y,v>0$. 
This $SL(2,R)$ subgroup is the ``$ax+b$ group'', namely the group of all 
linear transformations of the real line to itself that preserve orientation. 
With this framework, the operator $I_sG$ can be represented as a 
convolution operator 
$$I_sG = G* \psi_s
\tag 22$$
where convolution for left-invariant Haar measure on a locally compact group 
is defined by 
$$(f*g)(x) = \int_G f(y) g(y^{-1}x)\,dy\ .$$
Observe that $L_s(I_sG) =G$ for $s \ge 1$. 
The Riesz-Sobolev inequality and an extension of Young's inequality to 
non-unimodular groups provide good estimates for the fractional integral 
operator $I_s$. 
\enddemo

\subhead Riesz-Sobolev Inequality on $SL(2,R)/SO(2)$\endsubhead 
$$\int_M (f*g) (w) h(w)\,d\nu 
\le \int_M (f^* * g^*) (w) h^* (w)\, d\nu 
\tag 23$$ 
where $f,g$ and $h$ are non-negative measurable functions with $f^*,g^*$ and 
$h^*$ denoting their respective equimeasurable, geodesically decreasing 
rearrangements on $M\simeq SL(2,R)/SO(2)\break \simeq \IH^2$ and $d\nu$ is 
left-invariant Haar measure on $M$. 

\proclaim{Young's inequality}
Let $G$ be a locally compact group with left-invariant Haar measure denoted
by $m$. 
For $1\le p\le\infty$ 
$$\align
\|f*g\|_{L^p(G)} 
&\le \|f\|_{L^p(G)}\|\Delta^{-1/p'} g\|_{L^1 (G)}\\
\|f*g\|_{L^p(G)} 
&\le \|f\|_{L^1(G)}\| g\|_{L^p (G)}\\
\|f*g\|_{L^r(G)} 
&\le \|f\|_{L^p(G)}\|\Delta^{-1/p'} g\|_{L^q (G)} 
\endalign$$
where $\Delta$ denotes the modular function defined by $m(Ey) =\Delta (y)m(E)$,
$1/p + 1/p' =1$ and $1/r = 1/p + 1/q -1$. 
\endproclaim

\demo{Proof} 
Consider the form 
$$\align 
&\int_G h(x) (f*g)(x)\,dx = \int_{G\times G} h(x)f(y) g(y^{-1}x)\,dx\,dy\\
&\qquad 
= \int_{G\times G} h(x)f(xy) g(y^{-1})\,dx\,dy 
= \int_{G\times G} h(x)f(xy^{-1}) g(y)\Delta (y^{-1})\,dx\,dy\ .
\endalign$$
Then apply H\"older's inequality. 

It is natural here to consider fractional integration as a map from a space 
to its dual. 
The asymptotic behavior of $\psi_s$ combined with Young's inequality 
provide the necessary estimates to show that $I_s$ is a bounded map from 
$L^q (M)$ to $L^p(M)$ where $1/q + 1/p = 1$, $p=2+1/s$ and $q=2-1/(1+s)$. 
$$\alignat2
\psi_s (u) 
&\simeq \frac{\Gamma (s)\Gamma (s)}{4\pi\Gamma (2s)} u^{-s}
&&\qquad \text{as } u\to \infty\\
&\simeq - \frac1{4\pi} \ln u
&&\qquad \text{as } u\to 0\ .
\endalignat$$
Hence, any power of $\psi_s$ is locally integrable and using Young's 
inequality 
$$\|f*\psi_s\|_{L^p(M)} 
\le \|f\|_{L^q(M)} \|\Delta^{-1/p}\psi_s\|_{L^{p/2}(M)} 
= \|f\|_{L^q(M)} \|y^{1/p} \psi_s\|_{L^{p/2}(M)}\ .
\tag 24$$
The critical estimate is now reduced to the fact that $y^{s-1}(y+1)^{-2s}$ 
is integrable on $[0,\infty)$ for $s>0$. 
So the map $I_s$ is bounded from $L^q (M)$ to $L^p(M)$. 
The sharp constant for this estimate will be obtained using duality. 
\enddemo 

\proclaim{Theorem 7} 
For $s\ge 1$, $p=2+1/s$, $q= 2-1/(1+s)$ 
$$\gather
\|I_s G\|_{L^p(M)} \le A_p \|G\|_{L^q(M)} \tag 25\\
A_p = (2\pi)^{\frac2p -1} s^{-1-\frac2p}\ .
\endgather$$
This inequality is sharp and an extremal is given by $[1+d^2(z,i)]^{-1-s}$. 
For $F\in C^2(M)$ 
$$\Big[ \|F\|_{L^p(M)}\Big]^2 \le A_p \int_M F(L_sF)\,d\nu\ .
\tag 26$$ 
Here the extremal is $[1+d^2 (z,i)]^{-s}$. 
Because  $s\ge 1$, this latter result can be represented 
for $F\in C^2(M)\cap 
L^2(M)$ as 
$$\Big[ \|F\|_{L^p(M)}\Big]^2 \le A_p \biggl[ \int_M |DF|^2\,d\nu 
+ s(s-1) \int_M |F|^2\,d\nu\biggr] \ .
\tag 27$$
\endproclaim 

\demo{Proof} 
The plan of the argument is to use the Riesz-Sobolev inequality to show 
that an extremal function exists for (25) and hence by duality an extremal 
function exists for (26) which can be calculated using the 
Euler-Lagrange variational equation. 
To show the existence of an extremal for (25), it suffices to consider the 
functional 
$$\int_{M\times M} F(a) \psi_s [d^2 (z,w)] G(w)\,d\nu\,d\nu$$
for $F,G\ge0$ and $\|F\|_q = \|G\|_q =1$. 
By (24) this form is bounded above and by applying the Riesz-Sobolev 
inequality one can restrict attention to the case where $F$ and $G$ are 
geodesically radial decreasing functions. 
Then consider sequences of functions $\{F_n,G_n\}$ with $\|F_n\|_q = 
\|G_n\|_q =1$ so that 
$$\int_{M\times M} F_n(z) \psi_s [d^2 (z,w)] G_n(w)\,d\nu\,d\nu$$ 
converges to its maximum value. 
Since these functions are decreasing, one can use the Helly selection 
principle to choose  subsequences that converge almost everywhere to 
functions $F,G\in L^q(M)$. 
By Fatou's  lemma $\|F\|_q\le 1$, $\|G\|_q \le1$. 
Notice that $F_n(z) \le (4\pi u)^{-1/q}$, $G_n(z) \le (4\pi u)^{-1/q}$ 
using the radial variable $u=d^2 (z,i)$ since Haar measure restricted to the 
radial variable is $d\nu = 4\pi\,du$. 
Observe that 
$$[d(z,i)]^{-2/q} \psi_s [d^2 (z,w)] [d(w,i)]^{-2/q} 
\in L^1 (M\times M)\ .$$
Re-label the subsequences to have index $n$. 
By the dominated convergence theorem  
$$\int_{M\times M} F_n(z)\psi_s [d^2 (z,w)] G_n(w)\,d\nu\,d\nu 
\longrightarrow 
\int_{M\times M} F(z) \psi_s [d^2 (z,w)] G(w)\,d\nu\,d\nu$$ 
and so $\|F\|_{L^q(M)} = \|G\|_{L^q(M)} =1$ and $F,G$ must be extremal 
functions for (25). 
A somewhat similar argument is given in \cite3, page~40. 

{From} the Lemma above, one sees that if $G$ is an extremal for (25), then 
$F= I_sG$ is an extremal for (26).  
Moreover, if $G$ is radial decreasing, then $F$ will be radial decreasing 
since the convolution of two radial decreasing functions is radial decreasing. 
Hence, such an extremal $F$ must satisfy the Euler-Lagrange variational 
equation for (26): 
$$L_s F = \gamma F^{p-1}\quad ,\quad 
\|F\|_{L^p(M)} =1\ .$$
For $F$ a decreasing function of the radial variable $u= d^2 (z,i)$, one looks 
for solutions of the differential equation 
$$-\frac{d}{du} \left[ u(u+1) \frac{dF}{du}\right] + s(s-1) F 
= cF^{p-1}\ ,\qquad p= 2+1/s\ .
\tag 28$$ 
Note that if $F= I_sG$ for $G\in L^q(M)$ with $q=2-1/(1+s)$, then $F$ is 
bounded. 
Hence, there will be a unique solution to (28)  that is bounded and 
monotonically decreasing on $[0,\infty)$. 
This solution is 
$$F(u) = B(1+u)^{-s}$$ 
where the constant $B$ is determined by the condition that $\|F\|_q =1$. 
Now one can calculate the value of the sharp constant $A_p$. 
An extremal for (25) is obtained by 
$$L_sF = L_s (I_sG) =G\ .$$
This calculation completes the proof of Theorem~7. 
The argument developed here complements the result of Theorem~5.
Similar methods can also be applied for the case $0 < s < 1$ and will
be discussed in a more comprehensive treatment of Riesz potentials and
Sobolev embedding on hyperbolic space.
\enddemo

\head Acknowledgement\endhead

I would like to thank Tony Carbery for his invitation to participate in the 
Sussex Workshop on Fourier Analysis, and L.E.~Fraenkel, Nicola Garofalo, 
David Jerison and Eli Stein for useful remarks.


\Refs 

\ref\no 1
\by W. Beckner
\paper Sobolev inequalities, the Poisson semigroup and analysis
on the sphere $S^n$
\jour Proc. Nat. Acad. Sci. \vol 89 \yr 1992 \pages 4816--4819
\endref

\ref\no 2
\by W. Beckner
\paper Sharp Sobolev inequalities on the sphere and the
Moser-Trudinger inequality
\jour Ann. Math. \vol 138 \yr 1993 \pages 213--242
\endref

\ref\no 3
\by W. Beckner
\paper Geometric inequalities in Fourier analysis
\inbook Essays on Fourier Analysis in Honor of Elias M. Stein
\publ Princeton University Press
\yr 1995
\pages 36--68
\endref

\ref\no 4
\by W. Beckner
\paper Sharp inequalities and geometric manifolds
\jour J. Fourier Anal. Appl. \vol 3 \yr 1997 \pages 825--836
\endref

\ref\no 5
\by W. Beckner
\paper Geometric asymptotics and the logarithmic Sobolev inequality
\jour Forum Math. \vol 11 \yr 1999 \pages 105--137
\endref 

\ref\no 6 
\by H.J. Brascamp, E.H. Lieb and J.M. Luttinger
\paper A general rearrangement inequality for multiple integrals 
\jour J. Funct. Anal. 
\vol 17 \yr 1974 \pages 227--237 
\endref 

\ref\no 7 
\by G.B. Folland
\paper A fundamental solution for a subelliptic operator
\jour Bull. Amer. Math. Soc. \vol 79 \yr 1973 \pages 373--376
\endref

\ref\no 8
\by P.R. Garabedian
\book Partial differential equations \publ John Wiley \yr 1964
\endref

\ref\no 9 
\by V.V. Grushin 
\paper On a class of hypoelliptic operators 
\jour Math. Sbornik \vol 12 \yr 1970 \pages 458--475 
\endref

\ref\no 10 
\by S. Helgason
\book Differential geometry and symmetric spaces
\publ Academic Press
\yr 1962
\endref

\ref\no 11
\by S. Helgason
\paper Fundamental solutions of invariant differential operators on
symmetric spaces
\jour Amer. J. Math. \vol 86 \yr 1964 \pages 565--601
\endref 

\ref\no 12
\by D. Jerison and J.M. Lee
\paper Extremals for the Sobolev inequality on the Heisenberg group and
the CR Yamabe problem
\jour J. Amer. Math. Soc. \vol 1 \yr 1988 \pages 1--13
\endref

\ref\no 13
\by F. John
\paper The fundamental solution of linear elliptic differential equations
with analytic coefficients
\jour Comm. Pure Appl. Math. \vol 3 \yr 1950 \pages 273--304
\endref

\ref\no 14
\by S. Lang
\book $SL_2(R)$ \publ Addison-Wesley \yr 1975
\endref

\ref\no 15 
\by H. Lewy 
\paper An example of a smooth linear partial differential equation without 
solution 
\jour Ann. of Math. \vol 66 \yr 1957 \pages 155--158 
\endref 

\ref\no 16
\by A. Nagel and E.M. Stein
\book Lectures on pseudo-differental operators \publ Princeton University Press
\yr 1979
\endref

\ref\no 17
\by N. Varopoulos, L. Saloff-Coste and T. Coulhon
\book Analysis and geometry on groups \publ Cambridge University Press 
\yr 1992
\endref

\ref\no 18
\by N.J. Vilenkin
\book Special functions and the theory of group representations
\publ American Mathematical Society
\yr 1968
\endref

\ref\no 19
\by A. Weil
\book L'integration dans les groupes topologiques et ses applications 
\publ Hermann 
\yr 1966 
\endref

\ref\no 20
\by E.T. Whittaker and G.N. Watson 
\book A course of modern analysis
\publ Cambridge University Press 
\yr 1927 
\endref 

\endRefs
\enddocument